\documentclass{amsart}

\usepackage{euscript}   
\usepackage{amsfonts}
\usepackage{amsthm}
\usepackage{amsmath}
\usepackage{esvect}
\usepackage{amsfonts}
\usepackage{latexsym}
\usepackage{mathtools}
\usepackage{tikz}
\usepackage{graphicx}
\usepackage{tikz-cd}

\DeclareFontFamily{OT1}{pzc}{}
\DeclareFontShape{OT1}{pzc}{m}{it}{<-> s * [1.10] pzcmi7t}{}
\DeclareMathAlphabet{\mathpzc}{OT1}{pzc}{m}{it}  

\begin{document}

\theoremstyle{plain}
\newtheorem{CMS-lemma}{Lemma}[section]
\newtheorem{CMS-theo}[CMS-lemma]{Theorem}
\newtheorem{CMS-coro}[CMS-lemma]{Corollary}
\newtheorem{CMS-rema}[CMS-lemma]{Remark}
\newtheorem{CMS-propos}[CMS-lemma]{Proposition}
\newtheorem{CMS-conjecture}[CMS-lemma]{Conjecture}
\theoremstyle{plain}
\newtheorem{CMS-exa}[CMS-lemma]{Example}

\newcommand{\laba}[1]{\label{#1} \mbox{[   {#1}  ]}\ \ }

\def\bdem{\begin{proof}}
\def\edem{\end{proof}}
\def\bequ{\begin{equation}}
\def\eequ{\end{equation}}
\newcommand{\noi}{\noindent}
\newcommand{\papa}{\mbox{\(      H^{\infty}     \)}}
\newcommand{\chu}{\mathcal{M}(H^\infty)}
\newcommand{\rr}{\mbox{$    \longrightarrow   $}}
\newcommand{\clos}{\mbox{clos}}
\newcommand{\Om}{\Omega}
\renewcommand{\ker}{ \mbox{Ker}\, }
\newcommand{\range}{ \mbox{Ran}\, }
\newcommand{\la}{\langle}
\newcommand{\ra}{\rangle}
\newcommand{\disc}{ \mathbb{D} }
\newcommand{\re}{\mbox{Re}\,}
\renewcommand{\S}{U}  
\newcommand{\bolam}{\mathfrak{S}}
\newcommand{\Ker}{\text{Ker}\,}

\newcommand{\matru}[4]{
\left[     \begin{array}{cc}
          #1     &     #2   \\
          #3     &     #4
       \end{array}
\right]                }

\newcommand{\bN}{\mathbb{N}}
\newcommand{\bZ}{\mathbb{Z}}
\newcommand{\bT}{\mathbb{T}}
\newcommand{\bC}{\mathbb{C}}
\newcommand{\bR}{\mathbb{R}}

\newcommand{\cA}{\mathcal{A}}
\newcommand{\cH}{\mathcal{H}}
\newcommand{\cL}{\mathcal{L}}
\newcommand{\cJ}{\mathcal{J}}
\newcommand{\cN}{\mathcal{N}}
\newcommand{\cF}{\mathcal{F}}

\newcommand{\fM}{\mathfrak{M}}
\newcommand{\fD}{\mathfrak{D}}

\hyphenation{cha-rac-te-ri-za-tion}  \hyphenation{res-pec-ti-ve-ly}
\hyphenation{ve-ri-fy}  \hyphenation{na-tu-ral}
\hyphenation{co-ro-lla-ry}  \hyphenation{en-cou-ra-ge-ment}
\hyphenation{o-ther-wise}  \hyphenation{re-gu-la-ri-zed}
\hyphenation{sa-tis-fies}  \hyphenation{pa-ra-me-ters}
\hyphenation{subs-tan-tia-lly}  \hyphenation{des-com-po-si-tion}
\hyphenation{sur-pri-sing}  \hyphenation{pro-ducts}  \hyphenation{cons-truc-tion}
\hyphenation{theo-rem}  \hyphenation{o-pe-ra-tor}  \hyphenation{o-pe-ra-tors} \hyphenation{se-pa-ra-ted}
\hyphenation{eigen-va-lues}  \hyphenation{dia-go-na-li-zable} \hyphenation{fo-llows}
\hyphenation{in-ter-po-la-ting}  \hyphenation{a-ppli-ed}  \hyphenation{mul-ti-pli-ci-ties}
\hyphenation{in-va-riant}  \hyphenation{uni-lat-er-al} \hyphenation{a-ccep-ting}

\title{Frames of iterations and vector-valued model spaces}
\author{Carlos Cabrelli}
\author{Ursula Molter}
\address{{\em Carlos Cabrelli and Ursula Molter:} Departamento de Matem\'atica, Universidad de Buenos Aires, and Instituto de Matem\'atica ``Luis Santal\'o'' (IMAS-CONICET-UBA), 1428 Buenos Aires, Argentina}
\email{carlos.cabrelli@gmail.com, umolter@gmail.com}
\author{Daniel Su\'{a}rez}
\address{{\em Daniel Su\'arez:} Departamento de Matem\'atica, Universidad de Buenos Aires, and Instituto Argentino de Matem\'atica (IAM, CONICET)}
 \email{dsuarez@dm.uba.ar}
%
%

\begin{abstract}Let $T$ be a bounded operator on a Hilbert space $\cH$,
and $\cF = \{f_j , j \in J\}$ an at most countable set of vectors in $\cH$.
In this note we characterize the pairs
$(T, \cF)$ such that $\{T^nf : f \in \cF, n \in I\}$ form a frame of $\cH$, for
the cases of $I = \bN \cup \{0\}$ and $I = \bZ$.
\ \\
The characterization for unilateral iterations
gives a similarity with the compression of the shift acting on model spaces of
the Hardy space of analytic functions defined  on the unit disk
with values in $\ell^2(J)$. This generalizes recent work for iterations of a single function.
\ \\
In the case of bilateral iterations the characterization is by the bilateral shift acting
on doubly invariant subspaces of $L^2( T , \ell^2(J))$.
 \\
 Furthermore, we characterize the frames of iterations for vector-valued model operators when $J$ is finite
in terms of Toeplitz and multiplication operators in the unilateral and bilateral case, respectively. 
Finally we study the problem of finding the minimal number of orbits that produce a frame  in this context.
For the unilateral case we proved a formula for a lower bound. We conjecture that this lower bound is sharp. 
We give a new proof in the bilateral case, for which a formula is known.
\end{abstract}

\maketitle

\section{Introduction}

Let $\cH$ be a separable Hilbert space,  $T\in \cL (\cH)$ a bounded linear operator on $\cH$ and
$\{ x_j : \, j\in J\} \subset \cH$ be a subset of $\cH$, where the set of indexes $J$ is
countable (finite or countably infinite). A natural problem originated in dynamical  sampling asks for
a characterization of the operators and the vectors such that the orbits
$$
\cF= \{ T^n x_j : \ n\in \bN_0, \ j\in J \}
$$
is a frame, or even a Bessel sequence for $\cH$,  where $\bN_0$ denotes the set of nonnegative integers.

When $T$ is a normal operator and $\sharp J=1$ (the cardinal of $J$) the problem was solved in two papers
\cite{acmt} and \cite{accmp}, where the authors showed that $T$ is necessarily diagonalizable, the diagonal is
an interpolating sequence for the Hardy space  of the unit disk $H^2(\disc)$
(including finite interpolation if $\dim \cH$ is finite),
and the vector $x$ must satisfy some bound constraints depending on the eigenvalues of $T$.

In \cite{cmpp} the same problem is studied for an arbitrary finite set $J$. In this case the normal operator $T$ also must be diagonalizable, where the diagonal
is a union of at most $\sharp J$ interpolating sequences. A refinement of this characterization was obtained in \cite{CMS}, where the authors show that hyperbolic geometry appears naturally in the interplay
between the vectors $x_j\  (j\in J)$ and the sequence of eigenvalues.

In \cite{CHP} the authors consider the case of an arbitrary bounded operator $T$,
{\it not necessarily normal}, and   the iterations of a single vector $f\in \cH$
  i.e. $\#J = 1$.
They showed that  the orbit $\{T^nf : n \in \bN_0\}$ is a frame of $\cH$ if and only if
there exist a model space $K_{\theta}$ of $H^2(\disc)$, and a bounded invertible operator $V$,
from $K_{\theta}$ to $\cH$
such that $T$ is similar to $S_{\theta}$, i.e.  $T=VS_{\theta}V^{-1}.$ Here $S_{\theta}$ is the compression of the unilateral shift acting on $K_{\theta},$ and
$K_{\theta} = H^2\ominus \theta H^2$ with $\theta$ an inner function.
Furthermore, they show that there exists a function $\varphi \in K_{\theta}$ with $V\varphi=f$ such that $\{S^n_{\theta} \varphi: n \in \bN_0\}$ is a frame of $K_{\theta}.$

In other words, that theorem says that the orbit of a bounded operator in a Hilbert space $\cH$
forms a frame if and only if the orbit of the compression of the shift in
some model space isomorphic to $\cH$ forms a frame of the model space.
The authors also obtain similar results for the case of orbits  in $\bZ$,
and the space $L^2(\bT)$ for the characterization.

In  the present paper we show that these results can be extended  to the case of countable orbits.
For this we need to use model spaces in Hardy spaces with multiplicity.
We also generalized the case of $\bZ-$orbits, with the bilateral shift and reducing subspaces
of $L^2(\bT,\ell^2(J)).$

On the other hand, this characterization motivates the study of  frames of model spaces,
consisting of  orbits
of the compression of the shift.
Specifically, we consider the following problems.

For a given model space $\cN \subset H^2(\bT,\ell^2(J))$, characterize all the sets $\{\varphi_j: j\in J\} \subset \cN$
such that $\{S_{\cN}^n \varphi_j:j \in J, n \in \bN_0\}$ is a frame of $\cN.$ Here $S_{\cN}$
is the compression of the shift  to $\cN.$ We obtain necessary and sufficient conditions
in Theorem \ref{coronaframe}. We also address the same problem for $T$ invertible and $\bZ$-orbits, for which we find the necessary and sufficient conditions in Theorem \ref{2coronaframe},
and determine the minimum number of orbits needed to obtain a frame for $\cH$.

Finally,
given a bounded operator $T$ in a separable Hilbert space $\cH$,
which is the minimum cardinality of a set of functions $\cA$ in $\cH$
such that $\{T^n f: f \in \cA, n \in \bN_0\}$ is a frame of $\cH$?
Using the results mentioned above, this question can be studied in model spaces of some Hardy space. In this case we provide a conjecture (Conjecture~\ref{unifn}) together with an example supporting it.

The paper is organized as follows:

After establishing  the necessary notation and definitions, in Section~\ref{unilateral-sec} we
address the unilateral case. We first establish in Theorem~\ref{thebeginnig} a characterization of
when the forward iterations of a bounded operator on a set $\{f_j, j\in J\}$ in an arbitrary Hilbert space $\cH$ are a frame of $\cH$.

This motivates in Subsection~\ref{toep} the study of frames for the unilateral shift
of multiplicity $\sharp J$ acting on co-invariant  subspaces of $H^2(\mathbb{T},\ell^2(J))$.
In Subsection~\ref{minimum-subsec} we address the question of establishing the minimum cardinality
of unilateral orbits needed to form a frame.

Finally in Section~\ref{bilateral-sec} we look at doubly invariant subspaces, for which in Theorem~\ref{minimum-double} we establish the minimum number of orbits needed to obtain a frame for a doubly invariant subspace of of $(L^2)^m$.

\section{Unilateral orbits}\label{unilateral-sec}

Let $J$ be a countable set of indexes (i.e.: finite or infinite denumerable) and
denote by $H^2( \bT, \ell^2(J))$ the Hardy space of analytic functions from the
torus $\bT=\{ z\in \bC\!: |z|=1 \}$ into $\ell^2(J)$. That is, any $g\in H^2( \bT, \ell^2(J) )$
can be written as
\bequ\label{writes-as}
g(z)= \sum_{n\ge 0} a_n z^n,  \mbox{ where $a_n\in \ell^2(J)$ and }
\|g\|^2 = \sum_{n\ge 0} \|a_n\|_{\ell^2(J)}^2.
\eequ
Via this representation we can think of $\ell^2(J)$ as a closed subspace of
$H^2( \bT, \ell^2(J) )$.
Denote by $S$ the operator on $H^2( \bT, \ell^2(J) )$ of multiplication by $z$
(the unilateral shift  of multiplicity $\sharp J$). If $E\subset H^2( \bT, \ell^2(J) )$ is
a closed subspace and $P_E$ denotes the orthogonal projection onto $E$, the operator
 $$S_E:= P_E S|_E$$ is called the compression of $S$ to $E$.
We say that $E$ is coinvariant for $S$
if its orthogonal complement $E^\bot$ is $S$-invariant.
In this case, $E$ is called a model space and $S_E$ a model operator.

\begin{CMS-theo}\label{thebeginnig}
Let $\cH$ be a Hilbert space, $T\in\cL(\cH)$ and $f_j\in\cH$, where $j\in J$ is countable.
Then
\bequ\label{tn-de-fj}
\{ T^n f_j:  n\in \mathbb{N}_0 , \ j\in J\}
\eequ
is a frame for $\cH$ if and only if there is an $S$-coinvariant subspace
$K\subset H^2( \bT, \ell^2(J) )$  and an invertible operator
$V: K \to \cH$ such that
$$
V^{-1} T V = S_K \ \mbox{ and }\ V (P_K e_j) =  f_j   \ \mbox{ for all  $j\in J$},
$$
where $e_j \in \ell^2(J)$  is the standard basis.
\end{CMS-theo}
\bdem
Since $\{S^n e_j: \, n\in \mathbb{N}_0, \, j\in J \}$  is an ortonormal basis of
$H^2( \bT, \ell^2(J) )$,
we can define a linear transformation $W: H^2( \bT, \ell^2(J) ) \to \cH$ by setting
\bequ\label{w-del-shift}
W(S^n e_j) = T^n f_j ,  \  \ \forall n\in \mathbb{N}_0, \, j\in J .
\eequ
The operator $W$ is then bounded because  \eqref{tn-de-fj} is a Bessel sequence and it is onto because
this sequence is a frame. Therefore, letting $K:= (\ker W)^\bot$, we have that
$W|_{K}: K \to \cH$ is invertible.
In addition, it is clear from \eqref{w-del-shift} that $TW=WS$, implying that $\ker W$ is $S$-invariant,
which means that $K$ is coinvariant. Consider the
compression $S_K$ of $S$ to $K$. We claim that the following diagram commutes
\begin{equation*}
 \begin{tikzcd}
   \cH \arrow{r}{T} & \cH \arrow{d}{W|_K^{-1}} \\
    K  \arrow{u}{W|_K}\ \arrow{r}{S_K}& K
  \end{tikzcd}
\end{equation*}
First notice that $W|_K P_K (S^ne_j)=W(S^ne_j)= T^nf_j$. Hence, applying $T$ we get
$$ T W|_K P_K (S^ne_j)= T^{n+1} f_j = W|_K P_K (S^{n+1} e_j).
$$
Consequently,
\begin{align*}
W|_K^{-1} T W|_K P_K (S^ne_j) &=  P_K (S^{n+1} e_j) =  (P_K S)  (S^{n} e_j) \\
&=(P_K S)  (P_K S^{n} e_j) +   (P_K S)  (P_{K^\bot}S^{n} e_j) \\
&\stackrel{(1)}{=}  (P_K S)  (P_K S^{n} e_j),
\end{align*}
where (1) holds because $K^\bot$ is $S$-invariant. Since $P_K S^{n} e_j$
generate $K$, this proves our claim. Thus, the operator $V:= W|_K$ satisfies the theorem.
\edem
Now, we have the following proposition whose proof is immediate.

\begin{CMS-propos}\label{prop}
Let $\cH_0, \, \cH$ be Hilbert spaces, $T\in\cL(\cH)$ and $f_j\in\cH$, for $j\in J$ countable.
If $V:  \cH_0 \to \cH$ is an invertible operator,  then
$
\{ T^n f_j:  n\in \mathbb{N}_0 , \ j\in J\} \hspace{2mm} \mbox{is a frame for $\cH$}
$
if and only if $\displaystyle \, \{ (V^{-1} T V )^n (V^{-1}f_j):  n\in \mathbb{N}_0 , \ j\in J\}$
is a frame for $\cH_0$.
\end{CMS-propos}
Thus, combining  Proposition \ref{prop} with Theorem \ref{thebeginnig} we have the following corollary.

\begin{CMS-coro}\label{coro1}
Let $\cH$ be a Hilbert space, and $\{f_j\}_{j\in J}$ a set of vectors in $\cH$.The following assertions are equivalent.
\begin{enumerate}
\item $\{ T^n f_j:  n\in \mathbb{N}_0 , \ j\in J\} \hspace{2mm} \mbox{is a frame for $\cH$}$.
\item There exists an $S$-covariant space $K$ in $H^2( \bT, \ell^2(J))$, a collection of functions $\{\varphi_j\}_{j\in J}\subset K$ and, an isomorphism $V: K \to \cH$ such that

$\{ S_K^n \varphi_j:  n\in \mathbb{N}_0 , \ j\in J\} \hspace{2mm} \mbox{is a frame for $K$}$
with  $V^{-1} T V = S_K \ \mbox{ and }\ V (\varphi_j) =  f_j   \ \mbox{ for all  } j\in J$.
\end{enumerate}
\end{CMS-coro}

In short, what Corollary \ref{coro1} says is that every frame formed by iterations of a bounded operator
applied to a countable set of a Hilbert space can be
represented as a frame of iterations of a model operator.

\subsection{Toeplitz operators on model spaces}\label{toep}

In this section we will only consider the case where the set $J$ of Theorem \ref{thebeginnig} is finite.
By that theorem we can transfer the study of frames from an abstract context to the more familiar
scheme of model spaces. There we have many powerful tools from operator and function theory that we
can use to characterize Bessel sequences and frames of orbits. After doing so, we can go back to the original setting
by the invertible operator $V$ that implements the similarity in Theorem \ref{thebeginnig}.

One of the main tools at our disposal are Toeplitz operators with matricial symbols. In this section we shall
see that the orbits generated by a finite set $J$ form a Bessel sequence or a frame if an only if
a certain associated  Toeplitz operator is bounded or bounded above and below (respectively). Furthermore, we find
simple necessary and sufficient conditions for the respective boundedness of this operator.\\

For $m, k\in \bN$,  $L^2=L^2(\bT)$ and a matrix $G\in (L^2)^{m \times k}$ let
$$D_G =\{\vv{h}\in (H^2)^k  : \  G\vv{h}\in (L^2)^m \},
$$
where $\vv{h}$ is taken as a column vector.
Then for any linear submanifold $D\subset  D_G$ we can define a Toeplitz operator
$T_G: D\subset (H^2)^{k} \to (H^2)^{m}$ by the rule
$$
T_G \vv{f} = P_{(H^2)^{m}} (G \vv{h}),
$$
where $P_{(H^2)^{m}}$ is the orthogonal projection.
Observe that $D_G$ is dense in $(H^2)^k$, since it contains $(H^\infty)^k$.
If we write  $G= (g_{ij})$ it follows easily that $T_G$ extends boundedly to $(H^2)^k$ if and only if
all the scalar operators $T_{g_{ij}}$ extend boundedly to $H^2$, which is equivalent to $g_{ij}\in L^\infty$
for all $i,j$ (see \cite[Prop.$\,$10.1]{zhu}).
However, if $T_G$ is defined on a linear submanifold $D\subset D_G$,
the above condition could be unnecessary for extending $T_G$ continuously to the closure of $D$.
If $G= (g_{ij})\in (L^2)^{m\times k}$ write $G^\ast : = (\overline{g}_{ji})\in (L^2)^{k\times m}$
for the conjugate transpose  matrix. A simple calculation  shows that
$$
\la T_G \vv{h}, \vv{f} \ra = \la  \vv{h}, T_{G^\ast} \vv{f} \ra
\ \ \mbox{ for }\,  \vv{h}\in D_G , \, \vv{f}\in D_{G^\ast},
$$
meaning that $T^\ast_{G} =T_{G^\ast}$.
In particular, given a column matrix $\vv{g}=(g_1, \ldots , g_m)\in (H^2)^{m \times 1}$,
the Toeplitz operator $T_{\vv{g}^\ast} \vv{h}\in H^2$
for $\vv{h} \in (H^\infty)^m$ is densely defined in $(H^2)^m$.
The (unilateral) shift $S$ of multiplicity $m\in \mathbb{N}$ is the Toeplitz
operator with symbol $z I_{\mathbb{C}^m}$ acting on $(H^2)^m$.
That is, $S$ multiplies each coordinate of $h\in (H^2)^m$ by the variable $z$. It is clear that
$ST_\Theta = T_\Theta S$ when $\Theta\in (H^\infty)^{m\times m}$, implying that the
closure of $T_\Theta(H^2)^m$ is a closed $S$-invariant subspace.
The Beurling-Lax-Halmos Theorem
(for the respective scalar \cite{beur}, finite dimensional \cite{lax} and infinite dimensional
\cite{halm} case) says that  all the $S$-invariant subspaces are of this form and
they can be written in a  more precise and systematic form.

A matrix $\Theta\in (H^\infty)^{m\times m}$ is called rigid if the non-tangential limits
$\Theta (e^{it})$ are partial isometries on $\mathbb{C}^m$ with the same initial space for almost
every point $e^{it}\in \bT$. We include the null matrix as a particular case,
and we point out that when the initial space is all of $\mathbb{C}^m$ the matrix is called inner.
It follows that $T_\Theta (H^2)^m$ is closed and therefore an $S$-invariant subspace.
The B-L-H Theorem says that these are all.
Consequently, all the $S$-coinvariant subspaces (i.e.: $S^\ast$-invariant) have the form
$$
K_\Theta := (H^2)^m \ominus T_{\Theta} (H^2)^m
$$
for a rigid matrix $\Theta\in (H^\infty)^{m\times m}$.
As stated before in a more abstract setting, if $P_{K_\Theta}= I-T_{\Theta}T_{\Theta^\ast}$ denotes the orthogonal projection onto $K_\Theta$, the
operator obtained by compressing $S$ to this space:
$$
S_\Theta :=P_{K_\Theta} S \big|_{K_\Theta} : K_\Theta \to K_\Theta
$$
is  called a model operator, while $K_\Theta$ is a called a model space.
Since polynomials $\mathcal{P}^m$ are dense in $(H^2)^m$, the linear manifold
$$
(I-T_{\Theta}T_{\Theta^\ast})  (\mathcal{P}^m) \subset (H^\infty)^m \cap K_\Theta
\ \mbox{ is dense in $K_\Theta$}.
$$
Therefore, when $\vv{g}\in (H^2)^m$, the restriction of $T_{\vv{g}^\ast}$
to $(H^\infty)^m \cap K_\Theta$ is also densely defined in $K_\Theta$.
The next lemma characterizes the symbols $\vv{g}$ that make it bounded.

\begin{CMS-lemma}\label{Don}
$T_{\vv{g}^\ast} $ is bounded on $K_\Theta$ if and only if
\bequ\label{sarasony}
\vv{g}\in (H^\infty)^m + T_{\Theta} (H^2)^m .
\eequ
Moreover, if $\vv{g} \in T_{\Theta} (H^2)^m$ then $T_{\vv{g}^\ast}\equiv 0$ on $K_\Theta$.
\end{CMS-lemma}
\bdem
Let us momentarily denote the backward shift of multiplicity $m$ by $S^\ast_m$ .
It is easy to deduce the intertwining relation
$T_{\vv{g}^\ast} S^\ast_m = S^\ast_1 T_{\vv{g}^\ast}$
on  $(H^\infty)^m \cap K_\Theta$,
and if $T_{\vv{g}^\ast}|_{K_\Theta}$ is bounded, this extends to every $\vv{h}\in  K_\Theta$.
Now we can use a powerful tool of operator theory, the Sz.-Nagy and Foias
Commutant Lifting Theorem (see \cite[Thm.$\ $5]{d-m-p}), allowing us to lift the operator
$T_{\vv{g}^\ast}|_{K_\Theta}$ to a bounded operator on $(H^2)^m$ that also satisfies the
intertwining relation.
So, assuming boundedness, the mentioned theorem
tells us that there is a bounded operator $L: (H^2)^m \to H^2$ such that
$$
L|_{K_\Theta} = T_{\vv{g}^\ast}|_{K_\Theta}
\ \mbox{ and }\
L S^\ast_m  = S^\ast_1 L  \mbox{ on $(H^2)^m$}.
$$
Hence, the restriction of $L$ to each coordinate of $(H^2)^m$ is a bounded operator on $H^2$ that
commutes with
$S^\ast_1$, meaning that $L= T_{\vv{f}^\ast}$ for some $\vv{f}\in (H^\infty)^{m \times 1}$.
So, if we write $\vv{\mathsf{g}}= \vv{g}-\vv{f}$ then $T_{\vv{\mathsf{g}}^\ast} \equiv 0$ on $K_\Theta$,
and consequently when $\vv{h}\in K_\Theta$:
$$
0=(T_{\vv{\mathsf{g}}^\ast}\vv{h})(0) =
P_{H^2}   \big( \overline{\mathsf{g}}_1  h_1  + \cdots +  \overline{\mathsf{g}}_m   h_m \big) (0)
= \la \vv{h}, \vv{\mathsf{g}}\ra_{(H^2)^m}.
$$
This means that $\vv{\mathsf{g}}\in (K_\Theta)^\bot = T_{\Theta} (H^2)^m$,
from which \eqref{sarasony} follows immediately.

The reciprocal implication follows from the last assertion of the Lemma.
Therefore it only remains to show that if $\vv{g}= T_{\Theta} \vv{d}$, with
$\vv{d}\in (H^2)^m$, then $T_{\vv{g}^\ast}$ vanishes on $K_\Theta$.
For $f\in H^\infty$ and $\vv{h}\in(H^\infty)^m \cap K_\Theta$ we see that
\begin{align*}
\la T_{\vv{g}^\ast} \vv{h} , f\ra  &= \la \vv{h}, \vv{g} f \ra_{(H^2)^m}
= \la \vv{h}, T_{\Theta} (\vv{d} f) \ra_{(H^2)^m}
= \la T_{\Theta^\ast} \vv{h}, \vv{d} f \ra_{(H^2)^m}=0,
\end{align*}
where the last equality holds because $K_\Theta = \ker  T_{\Theta^\ast}$.
\edem

\vspace{2mm}\noi
Suppose that $\vv{g}\in K_\Theta \subset (H^2)^m$. Then for $\vv{h}\in (H^\infty)^m \cap K_\Theta$:
\begin{align*}
\la  \vv{h}, S_\Theta^n \vv{g} \ra &= \sum_{j=1}^m \la h_j, z^n g_j \ra
= \big\la \sum_{j=1}^m T_{\overline{g}_j}h_j  , z^n  \big\ra = \la T_{\vv{g}^\ast} \vv{h}  , z^n  \ra .
\end{align*}
Thus,
$\sum_{n \ge 0} |\la \vv{h}, S_\Theta^n \vv{g} \ra|^2 
=\| T_{\vv{g}^\ast}\vv{h}\|^2$. Furthermore, if $\vv{g}_1 , \ldots , \vv{g}_k\in K_\Theta$
are  column vectors and we consider the matrices
$$
G:= \big[\vv{g}_{1}, \ldots , \vv{g}_{k}\big]
\in  (H^2)^{m\times k}
\ \ \ \mbox{ and }\ \ \
T_{G^\ast} \vv{h} =
  \begin{bmatrix}
    T_{ \vv{g}_{1}^\ast }\vv{h}\\
    \vdots \\
    T_{ \vv{g}_{k}^\ast }\vv{h} \\
  \end{bmatrix}
\in  (H^2)^{k\times 1} ,
$$
it follows that
\bequ\label{bess2-frame}
\sum_{j= 1}^k \sum_{n \ge 0} |\la \vv{h}, S_\Theta^n \vv{g}_j \ra|^2
=\sum_{j= 1}^k  \| T_{\vv{g}_j^\ast}\vv{h}\|^2 = \|  T_{G^\ast} \vv{h}  \|^2 .
\eequ

\noi Hence, by \eqref{bess2-frame} the family
$\mathcal{F}:=\{ S_\Theta^n \vv{g}_j : \ n\in \mathbb{N}_0 , \ 1\le j\le k \}$
is a Bessel sequence in $K_\Theta$ if and only if  $T_{G^\ast} \big|_{K_\Theta}$ is bounded,
and is a frame if and only if it is  bounded above and below.

\begin{CMS-lemma}\label{tujess}
Let $F\in (H^\infty)^{m\times k}$ and $\Theta$ be a rigid matrix. Then
$$
T_{F}T_{F^\ast} \ge \delta^2 I_{(H^2)^m}\mbox{ on } K_\Theta \hspace{3mm}\mbox{ for some $\delta>0$}
$$
if and only if
$\displaystyle \   \|T_{F^\ast}h \|   +   \|T_{\Theta^\ast}h \|    \ge \varepsilon\|h\|$
for all $h\in (H^2)^m$ and some $\varepsilon>0$.
\end{CMS-lemma}
\bdem
The implication ($\Leftarrow$) follows by noticing that $\ker T_{\Theta^\ast}= K_\Theta$.
To prove ($\Rightarrow$) we can assume that $\|T_{F^\ast}\| \le 1$.
If $h\in (H^2)^m$ has  norm 1,
\begin{align*}
\delta \big[ 1 -   \|T_{\Theta^\ast}h \| \big]
&\le  \delta \| (I- T_{\Theta}T_{\Theta^\ast}) h \|
\ \le  \| T_{F^\ast}(I- T_{\Theta}T_{\Theta^\ast}) h \| \\*[1mm]
&\le  \  \| T_{F^\ast}h\|  +    \| T_{F^\ast}T_{\Theta}T_{\Theta^\ast} h \|
\ \le   \| T_{F^\ast}h\|  +    \| T_{\Theta^\ast} h \|,
\end{align*}
where the first inequality follows from $\|T_\Theta\| \le 1$ and the second inequality holds by hypothesis.
If $\|T_{\Theta^\ast}h \| \le  1/2$ this proves that
$\delta/2  \le \| T_{F^\ast}h\|  +    \| T_{\Theta^\ast} h \|$.
Otherwise $\|T_{\Theta^\ast}h \| >  1/2$ and there is nothing to prove.
\edem
\vspace{2mm}

\noi
The next Proposition requires a version of the Matrix Corona Theorem proved independently by Fuhrmann
\cite[Thm.$\,$3.1]{fuh} and Vasyunin (credited in \cite[Thm.$\,$3]{tolo}).

\begin{CMS-theo}\label{matrixcorona}
Let  $F_j\in (H^\infty)^{m\times m}$
for $1\le j\le q$. Then there are $B_j\in (H^\infty)^{m\times m}$ such that
$$\sum_{j=1}^q F_j(z)B_j(z) = I_{\mathbb{C}^m} \ \ \ \forall z\in \mathbb{D}.
$$
if and only if there is some $\eta>0$ such that
$\, \sum_{j=1}^q F_j(z) F^\ast_j(z)   \ge \eta^2 I_{\mathbb{C}^m}$,
$\ \ \forall z\in \mathbb{D}$.
\end{CMS-theo}

\begin{CMS-propos}\label{manyy}
Let $\Theta\in (H^\infty)^{m\times m}$ be a rigid matrix and $F\in (H^\infty)^{m\times k}$.
The following assertions are equivalent.
\begin{enumerate}
\item[{\em (1)}] $T_{ F } T_{ F ^\ast} \ge \delta^2 I_{(H^2)^m}$ on $K_\Theta$ for some $\delta>0$.
\item[{\em (2)}] $F (z)  F ^\ast(z) + \Theta(z) \Theta^\ast(z) \ge \eta^2 I_{\mathbb{C}^m}$
for some $\eta>0$ and all $z\in \mathbb{D}$.
\item[{\em (3)}] There exist $A\in (H^\infty)^{k \times m}$ and $B\in (H^\infty)^{m\times m}$ such that
$$
 F (z) A(z) + \Theta(z) B(z) = I_{\mathbb{C}^m}  \ \ \ \forall z\in \mathbb{D}.
$$
\end{enumerate}
\vspace{-2mm}
\end{CMS-propos}
\bdem
By Lemma \ref{tujess}, condition (1) says that there is some $\varepsilon >0$ such that
\bequ\label{tefeta}
T_{ F } T_{ F ^\ast} + T_{\Theta} T_{\Theta^\ast} \ge \varepsilon^2 I_{(H^2)^m}.
\eequ
For $\vv{c}\in \mathbb{C}^m$ consider $\vv{c}_{\!z} :=\vv{c}k_z \in (H^2)^m$,
where $k_z\in H^2$ is the normalized reproducing kernel of
$z\in \mathbb{D}$. Then $T_{ F ^\ast} \vv{c}_{\!z} =  F ^\ast(z) \vv{c} k_z$, from which
$$
\|T_{ F ^\ast} \vv{c}_{\!z} \|_{(H^2)^k}  
= \| F ^\ast(z) \vv{c} \|_{\mathbb{C}^k}
\ \mbox{ and }\
\|T_{\Theta^\ast} \vv{c}_{\!z} \|_{(H^2)^m}= \|\Theta^\ast(z) \vv{c} \|_{\mathbb{C}^m}.
$$
These equalities with \eqref{tefeta} give
$\displaystyle \| F ^\ast(z) \vv{c} \|^2_{\mathbb{C}^k} + \|\Theta^\ast(z) \vv{c} \|^2_{\mathbb{C}^m}
\ge \varepsilon^2 \|\vv{c} \|^2_{\mathbb{C}^m}$,
for all $\vv{c} \in\mathbb{C}^m$ and $z\in\disc$, which is (2) with $\eta = \varepsilon$.

Now suppose that (2) holds. If $k=m$ then (3) follows immediately from 
Theorem \ref{matrixcorona}.
Otherwise the idea is to use a very simple trick of linear algebra to reduce our problem to that theorem.
Write $k=mp +r$, where $p\in \mathbb{N}_0$,  $\,0\le r < m$, and decompose the matrix $F(z)$ into
$p$ square matrices $F_j(z)\in (H^\infty)^{m\times m}$ formed by each of the columns of $F(z)$, and
if $r>0$ include another matrix $F_{p+1}(z)$ formed by the last $r$ columns of $F(z)$ plus
columns of zeros until completing  a $(m\times m)$-matrix. Assume this  case to avoid
changing indexes. Then condition (2) translates into
$$
\sum_{j=1}^{p+1}F_j(z) F^\ast_j(z) + \Theta(z) \Theta^\ast(z) \ge \eta^2 I_{\mathbb{C}^m},
\ \ \forall z\in \mathbb{D}.$$
Thus, Theorem \ref{matrixcorona} says that there are matrices $A_j , B \in (H^\infty)^{m\times m}$ such that
$$
\sum_{j=1}^{p+1}F_j(z) A_j(z) + \Theta(z) B(z) = I_{\mathbb{C}^m},  \ \ \ \forall z\in \mathbb{D}.$$
It is clear then that if $A\in (H^\infty)^{k\times m}$ is formed by all the files of $A_1$ followed by
the files of $A_2$ and so on until the first $r$ files of $A_{p+1}$, then $A$ and $B$ satisfy  (3).

Finally, if (3) holds, passing to Toeplitz operators we get
$$
I_{(H^2)^m} = T_{ F  A + \Theta B }^\ast
= (T_{ F } T_A + T_{\Theta} T_B)^\ast =
T_{A^\ast} T_{ F ^\ast} + T_{B^\ast} T_{\Theta^\ast},$$
which implies that $T_{ F ^\ast}|_{K_\Theta}$ has a left inverse, and this is equivalent to (1).\\
\edem

\begin{CMS-theo}\label{coronaframe}
Let $\Theta\in (H^\infty)^{m\times m}$ be a rigid matrix and  $\vv{g}_j \in K_\Theta$ for $j=1, \ldots , k$.
Then the family
$$
\mathcal{F}=\{ S_\Theta^n \vv{g}_j: \, n\in \mathbb{N}_0 , \, 1\le j\le k \}$$
is a Bessel sequence in $K_\Theta$ if and only if
$\vv{g}_j= P_{K_\Theta} (\vv{f}_j)$, with $\vv{f}_j\in (H^\infty)^m$ for all $1\le j\le k$,
and a frame for $K_\Theta$  if and only if, in addition,
$F:= \big[\vv{f}_{1}, \ldots , \vv{f}_{k}\big]\in  (H^\infty)^{m\times k}$
satisfies any of the conditions of Proposition \ref{manyy}.
\end{CMS-theo}
\bdem
Consider the matrix $\displaystyle G:= \big[\vv{g}_{1}, \ldots , \vv{g}_{k}\big] \in  (H^2)^{m\times k}$,
where each $\vv{g}_j$ is a column vector.
By the comment that follows \eqref{bess2-frame}, 
$\mathcal{F}$ is a Bessel sequence (or a frame) for $K_\Theta$
if and only if $T_{G^\ast} \big|_{K_\Theta}$ is bounded above (respectively above and below).
Since $\vv{g}_{j}\in K_\Theta$, Lemma \ref{Don} says that upper boundedness is equivalent to
the existence of $\vv{f}_j\in (H^\infty)^m$ such that $\vv{g}_j = P_{K_\Theta}\vv{f}_j$ for $1\le j\le k$.
So, writing $F= \big[\vv{f}_{1}, \ldots , \vv{f}_{k}\big] \in (H^\infty)^{m\times k}$, a new use of
Lemma \ref{Don} gives
$$
T_{G^\ast} \big|_{K_\Theta} = T_{F^\ast} \big|_{K_\Theta}.$$
Thus, $T_{G^\ast} \big|_{K_\Theta}$ is bounded below if and only if
$T_{F^\ast} \big|_{K_\Theta}$ satisfies (1) of Propos.$\ $\ref{manyy}.
\edem

\noi
Among the conditions of Proposition \ref{manyy} it is apparent that (2) is the easiest to check.
This fact will be reflected in the next two examples, where it is used to characterize frames
of orbits in accordance with Theorem \ref{coronaframe}.

\begin{CMS-exa}\label{exx1}
{\em
Suppose that $u\in H^2$ is inner or the null function and let $f_j \in H^\infty$ for $j=1,\ldots, k$,
be such that
$$
|f_1(z)|^2 + \cdots + |f_k(z)|^2 + |u(z)|^2 \ge \eta^2
$$
for some $\eta>0$ and all $z\in \mathbb{D}$, which is  (2) of Proposition \ref{manyy}.
Then the $S_u$ orbits of
$$g_j:= (I-T_u T_{\overline{u}}) f_j \in K_u \ \mbox{ for }\ j=1, \ldots k
$$
form a frame for $K_u$. By the last theorem any frame of at most $k$ orbits
(because of possible repetitions) has this form. Repetitions of the $g_j's$ can be avoided
by setting the additional condition $f_i-f_j\not\in u H^\infty$ for $1\le i < j \le m$.
Furthermore, if the set of functions $\{ f_j\!: 1\le j\le k \}$ that satisfies the above inequality is
minimal then so is $\{ g_j\!: 1\le j\le k \}$, in the sense that if we extract any of the respective orbits,
the remaining orbits no longer form a frame.

When $u\equiv 0$ the corresponding operator is the simple shift $S$ and a single
orbit is a frame if and only if it is generated by an invertible function $f\in \papa$.
Observe the contrast with the characterization of cyclic vectors for $S$ given by Beurling's theorem,
where the orbit  is generated by an outer function $g\in H^2$ (see \cite[Cor.$\ $7.3]{gar}).
}
\end{CMS-exa}

\begin{CMS-exa}\label{exx2}
{\em
For $b_1, b_2 \in H^\infty$ inner functions such that
$\inf_\disc |b_1|^2 + |b_2|^2 \ge \delta^2$ consider the inner matrix
$$
\Theta  =
  \begin{bmatrix}
    b_1 & 0\\
    0   & b_2 \\
  \end{bmatrix}
\in  (H^\infty)^{2\times 2}
$$
and the model operator $S_\Theta: K_\Theta \subset (H^2)^2 \to K_\Theta$.
If $F= {1 \choose 1} \in (H^\infty)^{2 \times 1}$ then
\begin{align*}
\left\| F^\ast(z) {c_1 \choose c_2} \right\|^2  + \left\| \Theta^\ast(z) {c_1 \choose c_2} \right\|^2
&=  |c_1 + c_2 |^2  + | b_1(z)c_1|^2   +  | b_2(z)c_2|^2 \\
&\ge  (|c_1| - |c_2 |)^2  + \delta^2 \min\{ | c_1|^2   , | c_2|^2 \} \\
&\ge  \frac{\delta^2}{16} \, ( | c_1|^2  +  | c_2|^2 )
\end{align*}
for all $(c_1,c_2)\in \mathbb{C}^2$ and $z\in\disc$, where the last inequality follows easily
by considering the two cases $2|c_2| \le |c_1|$ and $2|c_2| \ge |c_1|  \ge |c_2|$.
Since the above inequality is (2) of Proposition \ref{manyy},
then Theorem \ref{coronaframe} says that the $S_\Theta$-orbit generated by
$$
\vv{g}:= (I-T_\Theta T_{\Theta^\ast}){1 \choose 1} = {1-b_1 \overline{b_1(0)} \choose 1-b_2 \overline{b_2(0)}} \in K_\Theta
$$
is a frame for $K_\Theta$. Therefore $K_\Theta$ admits a frame of a single $S_\Theta$-orbit, which according
to Theorem \ref{thebeginnig} means that $S_\Theta$ is similar to a scalar model operator.
We can imitate the proof of that theorem to find the model space and the similarity.
Indeed, let $W: H^2 \to K_\Theta$ be the linear operator defined by
\vspace{-1mm}
$$
W(z^n) = S^n_\Theta \vv{g}= P_{K_\Theta} S^n {1\choose 1} .
$$
\vspace{-2mm}
Hence,
$\displaystyle W(h) =  P_{K_\Theta}  {h\choose h}$ for every $h\in H^2$, and consequently
$$
\Ker W = \Big\{ h\in H^2 : \, {h \choose h} \in T_\Theta H^2 \Big\}
= \{ h\in H^2 : \, h \in b_1 H^2 \cap b_2 H^2 \} = b_1 b_2 H^2 .
$$
So, $(\Ker W)^\bot = K_{b_1 b_2}$, the operator $V:= W|_{K_{b_1 b_2}} : K_{b_1 b_2} \to K_\Theta$ is invertible,
$$
V^{-1} S_\Theta V = S_{b_1b_2} \ \ \mbox{ and }\ \
V P_{K_{b_1b_2}}1 = \vv{g}.
$$
}
\end{CMS-exa}

\subsection{Minimum cardinality of unilateral orbits}\label{minimum-subsec}

\noi
{\bf Definition.}
Let $T: \mathcal{H}  \to \mathcal{H}$ be a bounded operator on a
separable Hilbert space $\mathcal{H}$.
We define the (unilateral) frame number $\vartheta_+(T)$ of $T$
as the minimum
$m\in \mathbb{N}\cup \{ \infty \}$ such that there exists a frame formed by the union of $m$ orbits.
If there is no $m$ with this property we say that $T$ has no unilateral frame number.
If $T$ is invertible, its (bilateral) frame number $\vartheta_-(T)$ is defined  analogously
with the bilateral orbits $\{T^n g: \, n\in \mathbb{Z} \}$, for $g\in \mathcal{H}$.
\vspace{3.5mm}

\noi Observe that Theorem \ref{thebeginnig} implies that  $T$ has frame number $m$ if and only if
it is similar to a model operator $S_\Theta$ acting on $K_\Theta \subset H^2(\bT, \ell^2(J))$,
where $m=\sharp J$ is the minimum with this property.
Examples \ref{exx1} and \ref{exx2} show two operators of frame number 1.
A simple example will show that a rank 1 operator could have infinite frame number or no number at all.

\begin{CMS-exa}\label{exx3}
{\em
Every rank 1 operator on the separable Hilbert space $\cH$ has the form
$Tx= \la x,g\ra f$ for all $x\in \cH$, where $f, g \in \mathcal{H}\setminus \{0\}$. Thus,
$$
T^0x=x \ \mbox{ and }\  T^nx = \la x, g\ra \, \la f, g\ra^{n-1} f,\/ \ \mbox{ for $n\in\mathbb{N}$. }
$$
Therefore
$$  \sum_{n\ge 0} |\la h, T^n x\ra|^2 = |\la h,  x\ra|^2 +
|\la x,  g\ra|^2  \, |\la h,  f\ra|^2  \,          \sum_{n\ge 1} | \la f, g\ra|^{2(n-1)}.
$$
%
If $| \la f, g\ra| \ge 1$, the series diverges, so when $x\not\perp g$ the orbit of $x$ is not even
a Bessel sequence (take $h=f$). When $x\bot g$, taking $h=g$ we see that $\la g , T^n x\ra =0$ for all $n$ and there
is no frame formed by orbits.
Finally, if $| \la f, g\ra| < 1$ and $\{  x_j :  j\in J\}$ is a frame, then so is
$\{ T^n x_j : n\in \mathbb{N}_0, \, j\in J\}$. That is,  $\vartheta_+(T)\le \sharp J$.
}
\end{CMS-exa}

\begin{CMS-rema}
{\em
From now on we write $\bolam:= \{\vv{c}\in \bC^m: \|\vv{c}\|=1\}$ for the unit sphere in $\bC^m$.
For $m\in \mathbb{N}$ let $\Theta\in (H^\infty)^{m\times m}$ be a rigid matrix
and consider the model operator
$S_\Theta : K_\Theta \to K_\Theta$. It is clear that if $\Theta \equiv 0$ then $S_\Theta$
actually is the shift of multiplicity $m$ and  $\vartheta_+(S_\Theta)= m$. This immediate observation
has an interesting generalization: if there is a sequence $z_n\in\disc$  such that
$\Theta (z_n) \to 0$ 
then $\vartheta_+(S_\Theta)= m$. Indeed, if $F\in (H^\infty)^{m\times k}$ satisfies
(2) of Propos.$\ $\ref{manyy} then
$$
\inf_{\vv{c}\in \bolam}\|   F ^\ast(z_n)  \vv{c} \|^2_{\mathbb{C}^m} \ge (\eta/2)^2 
\ \mbox{ for $n$ large enough.}
$$
This is not possible if $k<m$. This shows that in Example \ref{exx2} the hypothesis
$\inf_\disc |b_1|^2 + |b_2|^2 \ge \delta^2$ is necessary for the frame number to be 1.
This remark leads directly to Conjecture \ref{unifn} below.}
\end{CMS-rema}

\noi
If $T\in \cL(\cH)$ admits a frame formed by finitely many orbits, there is some $m\in\bN$ such that
$T$ is similar to a model operator $S_\Theta$ acting on $K_\Theta \subset  (H^2)^m$ for
some rigid matrix $\Theta \in (H^\infty)^{m \times m}$. In this case, Theorem \ref{coronaframe} tells us that
$\vartheta_+(T)= \vartheta_+(S_\Theta)=k$ is the minimum number of vectors $\vv{f}_j \in  (H^\infty)^{m}$,
with $1\le  j\le k$, such that if $F= [\vv{f}_1 \ldots \vv{f}_k] \in (H^\infty)^{m \times k}$, there
exists some $\eta >0$ satisfying (2) of Propos.$\ $\ref{manyy}.
If $\Theta \equiv 0$ (i.e.: $S_\Theta$ is the shift of multiplicity $m$) then
$F\in (H^\infty)^{m \times m}$ must be invertible. For  instance, $F$ is the identity in the proof
of Theorem \ref{thebeginnig}.

The maximal ideal space of $\papa$ is
\vspace{-2.5mm}
$$\chu = \{\varphi: \papa \to \bC\, \text{ linear, continuous, multiplicative and $\ne 0$}\}
$$
with the weak$^{\ast}$ topology, that makes it a compact Hausdorff space.
Via evaluations the disk can be viewed as an open subset of $\chu$ and thanks to the corona theorem of
Carleson \cite{cacor} we know that it is dense. Therefore $\chu$ is a compactification of $\disc$.
The Gelfand transform of $f\in \papa$ consists of looking at $f$ in the bidual space of $\papa$ and
then restricting its domain to $\chu$. Thus, the Gelfand transform provides a continuous
extension of $f$ to the whole $\chu$.
In particular,  the Gelfand transform induces a continuous extension of
any $\Theta \in (H^\infty)^{m\times m}$ as a map from $\chu$ into $\bC^{m\times m}$,
that we also denote by $\Theta$.

\begin{CMS-propos}\label{postt}
Let $0\ne K_\Theta\subset  (H^2)^m$ be the $S$-coinvariant subspace associated to the rigid matrix
$\Theta \in (H^\infty)^{m\times m}$.
Then
{\em
$$  
\vartheta_+ (S_\Theta) \ge  \sup_{x\in \chu }  \big(\! \dim \ker \Theta(x)\big) .
$$  
}
\end{CMS-propos}
\bdem 
We are looking for the minimum integer $k\in \bN$ such that there are $F\in (H^\infty)^{m\times k}$
and $\delta >0$ satisfying
$$
 \| F^\ast (z) \vv{c} \|_{\bC^k}^2  + \| \Theta^\ast (z) \vv{c} \|_{\bC^m}^2 \ge \delta
$$
for all $z\in\disc$ and $\vv{c}\in \bolam$, where $F^\ast(z): \bC^m \to \bC^k$.
By the density of $\disc$ in $\chu$ and the compactness of the latter, 
this is equivalent to
\bequ\label{maxcoro}
 \| F^\ast (x) \vv{c} \|_{\bC^k}^2  + \| \Theta^\ast (x) \vv{c} \|_{\bC^m}^2 >0,
 \ \ \forall x\in \chu     ,   \ \forall \vv{c}\in\bolam .
\eequ
By the dimension theorem, the kernels of $\Theta(x)$ and $\Theta^\ast(x)$ have the same dimension.
Let $p$ be the right number in \eqref{framechu} and $x\in \chu$ be such that
$\dim \ker \Theta^\ast(x)=p$. If $k$ and $F$ are such that \eqref{maxcoro} holds and we restrict
$$
F^\ast(x)|_{\text{Ker} \,\Theta^\ast(x)}: \ker \Theta^\ast(x)\approx \bC^p \to \bC^k,
$$
the dimension theorem tells us that
$$p= \dim \big(\ker F^\ast(x)|_{\text{Ker} \,\Theta^\ast(x)}\big) +
\dim \big(\range F^\ast(x)|_{\text{Ker} \, \Theta^\ast(x)}\big),
$$
where the last summand is $\le k$. If $k<p$ the first summand must be nonzero, which implies that
$\ker F^\ast(x)\cap \ker \Theta^\ast(x)\neq 0$. Since this contradicts \eqref{maxcoro}, we get $k\ge p$.
\edem

The only reason for having a single inequality in Proposition \ref{postt} is that we were unable
to prove the other one. Below we conjecture the equality, which provides a formula for the unilateral
frame number. Furthermore, we show that it is enough to prove the conjecture
when $\Theta$ is an inner matrix, and then we give an example that supports the conjecture.

\begin{CMS-conjecture}[The unilateral frame number]\label{unifn}
Let $0\ne K_\Theta\subset  (H^2)^m$ be the $S$-coinvariant subspace associated to the rigid matrix
$\Theta \in (H^\infty)^{m\times m}$.
Then
{\em
\bequ\label{framechu}
\vartheta_+ (S_\Theta) =  \sup_{x\in \chu }  \big(\! \dim \ker \Theta(x)\big) .
\eequ
}
\end{CMS-conjecture}

\begin{CMS-rema}
{\em
We shall see that it is enough to prove \eqref{framechu} 
when $\Theta$ is inner.
Only the inequality ($\le$) needs to be proved.
Denote by $E\subset \bC^m$ the common initial space of the partial isometries
$\Theta(e^{it})$ for almost every $e^{it}\in \partial\disc$. If $\Theta$ is not an inner matrix,
$m_1:=\dim E <m$. Then there is an unitary matrix $V\in \bC^{m \times m}$ mapping
$\bC^{m_1} \times \{0\}^{m-m_1}$  onto $E$  and $\{0\}^{m_1} \times \bC^{m-m_1}$ onto $E^\bot$.
Consequently,
$$
\Theta V  =
  \begin{bmatrix}
     \Theta_1 & \hspace{0mm}     0 \\
  \!\!  0           & \hspace{0mm}  0\\
  \end{bmatrix}
 \in (\papa)^{m\times m},
  \text{ where }
  \Theta_1\in (\papa)^{m_1\times m_1}
  \text{ is inner.}
$$
Therefore, if there is $F_1\in (\papa)^{m_1\times k_1}$ (with $k_1 \le m_1$)
such that $F_1, \, \Theta_1$ satisfy \eqref{maxcoro} with $k_1$ and $m_1$ instead of $k$ and $m$, the pair
$$
F  =
  \begin{bmatrix}
     F_1 & \hspace{-3mm}     0 \\
  \!\!  0           & 
    I_{m-m_1}\\
  \end{bmatrix}
 \in (\papa)^{m\times (k_1+m-m_1)} \,\mbox{ and $\Theta V$ satisfy \eqref{maxcoro},}
$$
and the same holds for the pair $F$  and $\Theta$.}
\end{CMS-rema}

\noi
By the above reduction, to prove the conjecture we can assume that $\Theta$ is an inner matrix,
and consequently its determinant  is an inner function. We give below an example where, in addition,
the conjecture is true.

\begin{CMS-exa}
{\em The conjecture holds if $\Theta \in (H^\infty)^{m\times m}$ is such that $u:= \det \Theta$ is
an interpolating Blaschke product.

Indeed, consider $Z_u := \{ x\in \chu\!: \, u(x)=0  \}$.
If $p$ is the right number in \eqref{framechu}, for each $x\in Z_u$ there is
a constant matrix $F_x\in \bC^{m\times p}$ such that
$\| F_x^\ast \vv{c}\|_{\bC^p} \ge 1$ when $\vv{c}\in\ker \Theta^\ast (x)\cap \bolam$.
Thus, by the compactness of $\bolam$, for each fixed $x\in Z_u$ we have
$$
\beta_x := \inf_{\vv{c}\in \bolam}
\big( \| F_x^\ast \vv{c}\|_{\bC^p}^2 + \|\Theta^\ast(x) \vv{c}\|^2_{\bC^m} \big) >0,
$$
from which the following set is an open neighborhood of $x$ in $\chu$:
$$
V_x= \{ y\in\chu\!: \, \inf_{\vv{c}\in \bolam}
\big( \| F_x^\ast \vv{c}\|_{\bC^p}^2 + \|\Theta^\ast(y) \vv{c}\|^2_{\bC^m} \big) > \beta_x/2
\}
$$
Since $Z_u$ is compact, there is a finite covering by these sets $V_1, \ldots , V_N$,
and the associated matrices $F_j \in \bC^{m\times p}$, with $1\le j\le N$, such that for some $\alpha >0$,
\bequ\label{conW}
\inf_{\vv{c}\in \bolam}  \big( \| F_j^\ast \vv{c}\|_{\bC^p}^2
+ \|\Theta^\ast(y) \vv{c}\|^2_{\bC^m} \big) \ge \alpha \ \mbox{ when  $y\in  V_j$.}
\eequ
Furthermore, since  $Z_u$ is totally disconnected (see \cite[p.$\,$395]{gar}), by a well-known result of
dimension theory (see \cite[p.$\,$87]{nag}) we can assume that the $V_j's$ satisfy the extra condition
$$
(V_j\cap Z_u) \cap (V_k\cap Z_u) = \emptyset \ \mbox{ if } j\ne k .
$$
This means that every $y\in Z_u$ belongs to only one of the $V_j's$.
In particular, when $z_n\in Z_u\cap \disc$ let $V_j$ be the corresponding set such that $z_n\in V_j$.
Since the sequence $\{ z_n\}$ is interpolating, there is $F\in (\papa)^{m\times p}$
such that $F(z_n)= F_j$. By the density of $Z_u\cap \disc$ in $Z_u$ (see \cite[p.$\,$395]{gar}), for every $x\in Z_u$:
$$
F(x)= F_j,  \mbox{ where $V_j$ is the unique set that contains $x$.} 
$$
Consequently, \eqref{conW} says that
\begin{align*}
\inf_{\vv{c}\in \bolam} \big( \| F(y)^\ast \vv{c}\|_{\bC^p}^2 + \|\Theta^\ast(y) \vv{c}\|^2_{\bC^m} \big)
&\ge \alpha \ \mbox{ when  $y\in  \bigcup_{j=1}^N V_j$,}
\end{align*}
and since this  is an open neighborhood of $Z_u$, we also have
$\inf_{\vv{c}\in \bolam} \|\Theta^\ast(x) \vv{c}\|^2_{\bC^m} >0$ when $x\not\in  \bigcup_{j=1}^N V_j$.
That is, \eqref{maxcoro} holds for $F\in (\papa)^{m\times p}$, giving the inequality that is missing  in the
conjecture.}
\end{CMS-exa}

In the above example, $u$ could be a finite Blaschke product that satisfies finite interpolation
(i.e.: with single zeros). Moreover, it is clear that the argument in the Example can be easily adapted
to a finite Blaschke product times an interpolating Blaschke product.

\section{Bilateral orbits}\label{bilateral-sec}

Let $J$ be a countable set of indexes and let $\bZ$ be the set of integers.
Denote by $L^2( \bT, \ell^2(J) )$ the space of $L^2$ functions from $\bT$ into $\ell^2(J)$.
That is, any $g\in L^2( \bT, \ell^2(J) )$ writes as \eqref{writes-as},
where the sum runs over $\bZ$ and $z\in\bT$.
The bilateral shift $\S$ of multiplicity $\sharp J$ is the operator of multiplication by the variable $z$
on $L^2( \bT, \ell^2(J) )$. It is clear that $\S^\ast$ consists of multiplication by
$z^{-1}=\overline{z}$, so $\S$ is unitary.
A subspace $E\subset L^2( \bT, \ell^2(J) )$ is doubly invariant
if $E$ and $E^\bot$ are $\S$-invariant, or equivalently, $E$ is both $\S$ and $\S^\ast$ invariant.
A proof analogous to that of Theorem \ref{thebeginnig} gives

\begin{CMS-theo}             \label{genedoubly}
Let $\cH$ be a Hilbert space, $T\in\cL(\cH)$ be an invertible operator and $f_j\in\cH$,
where $j\in J$ is countable.
Then $\{ T^n f_j:  n\in \bZ , \ j\in J\}$ is a frame for $\cH$ if and only if there is
a doubly invariant subspace $K\subset L^2( \bT, \ell^2(J) )$  and an invertible operator $V: K \to \cH$
such that
$$
V^{-1} T V = \S|_K \ \mbox{ and }\ V (P_K e_j) =  f_j   \ \mbox{ for all  $j\in J$},
$$
where $e_j \in \ell^2(J)$ is the standard basis.
\end{CMS-theo}

\noi
The theorem relates general frames of bilateral orbits to doubly invariant subspaces of
$L^2( \bT, \ell^2(J) )$. The characterization of these spaces were obtained by Helson and Lowdenslager \cite{HL61} and
Srinivasan \cite{Sri64}.
We follow the treatment of Lecture VI in Helson's book \cite{hel}.
\vspace{1.6mm}

\noi
{\bf Definition.} A range function is a map defined almost everywhere
$$
\cJ : \bT \to \{ \mbox{closed subspaces of $\mathcal{H}$} \}$$
that is (weakly) measurable in the sense that if $P_{\cJ(z)}$ is the orthogonal projection on
$\cJ(z)$, then $z\mapsto\la P_{\cJ(z)}a, b\ra$ is measurable for all $a, b\in\mathcal{H}$.
\vspace{3.5mm}

\noi
In \cite[Thm.$\,$8]{hel} it is showed that a closed subspace $\fM \subset L^2(\bT, \mathcal{H})$ is
doubly invariant if and only if there is a (unique) measurable range function $\cJ$
such that
$$\fM  = \fM_\cJ:=
\{ f\in L^2(\bT,\mathcal{H}): \   f(z) \in \cJ(z),  \mbox{ a.e.}\ z\in \bT \}.
$$
So, the orthogonal projection $P_\fM : L^2(\bT,\mathcal{H}) \to \fM$ is given by
$$
(P_\fM \vv{f})(z)  =    P_{\cJ(z)} \vv{f}(z) \ \mbox{ for a.e.$\, z\in \bT$}.
$$
When $\mathcal{H} = \ell^2(J)$, where $\sharp J=m\in \mathbb{N}$,  we identify 
$L^2(\bT,\mathcal{H}) = L^2(\bT,\mathbb{C}^m) = L^2(\bT,\mathbb{C} )^m = (L^2)^m$,
where the last equality is just notation.
Therefore for $\mbox{a.e.}\, z\in \bT$,  the orthogonal projection
$P_{\cJ(z)}: \mathbb{C}^m \to \cJ(z)$ is represented by
a multiplication operator $M_{\sigma(z)}$ whose symbol is pointwise the matrix of a
projection in $\mathbb{C}^m$.
That is, there is a {\bf unique} matrix $\sigma\in (L^\infty)^{m\times m}$
such that $M_\sigma = P_{\fM}: (L^2)^m \to \fM$.
So, $P_{\cJ(z)}= M_{\sigma(z)}$ for a.e.$\,z\in \bT$.
Conversely, if $\sigma\in (L^\infty)^{m\times m}$ is any self-adjoint and idempotent matrix then
$\fM:= M_\sigma (L^2)^m$ is a doubly invariant subspace.
It follows immediately that $M_\sigma (L^\infty)^m =  (L^\infty)^m \cap \fM$ is dense in $\fM$.

\subsection{Multiplication operators on doubly invariant subspaces}

Most of the results in this subsection are completely analogous to those for Toeplitz operators
on model subspaces, except that the multiplication operator is not followed by the analytic projection.
This greatly simplifies many arguments.
\vspace{1mm}

For  $m, k\in \bN$,  $L^2=L^2(\bT)$ and $G\in (L^2)^{m \times k}$ define the multiplication operator
$$
M_G: (L^\infty)^{k} \subset (L^2)^{k} \to (L^2)^{m} \ \mbox{ as }\ M_G \vv{f} = G \vv{h}.
$$
It is densely defined on $(L^2)^{k}$ and if we write $G= (g_{ij})$, it is immediate that
$M_G$ is bounded on $(L^2)^k$ if and only if the scalar operators $M_{g_{ij}}$
are bounded on $L^2$, which is equivalent to $g_{ij}\in L^\infty$ for all $i,j$.
Also, $M^\ast_{G} =M_{G^\ast}$.

Let $\fM\subset (L^2)^m$ be a doubly invariant subspace and $\vv{g}\in (L^2)^m$.
Then $M_{\vv{g}^\ast}$     
is densely defined on $\fM$, and for $\vv{h}\in (L^\infty)^m \cap \fM$:
\bequ\label{bibess0}
\la  \vv{h}, \S^n \vv{g} \ra
=  \sum_{i=1}^m \la h_i, z^n g_i \ra   
= \big\la \sum_{i=1}^m  M_{\overline{g}_i}  h_i  , z^n  \big\ra = \la M_{\vv{g}^\ast} \vv{h}  , z^n  \ra .
\eequ

\noi
Thus, if $\vv{g}_1 , \ldots , \vv{g}_k\in (L^2)^{m}$ are the column vectors of  $G\in  (L^2)^{m\times k}$,
\bequ\label{bibess2-frame}
\sum_{j= 1}^k \sum_{n \in \bZ} |\la \vv{h}, \S^n \vv{g}_j \ra|^2
=\sum_{j= 1}^k  \| M_{\vv{g}_j^\ast}\vv{h}\|^2 = \|  M_{G^\ast} \vv{h}  \|^2 .
\eequ

\begin{CMS-lemma}\label{biDon}
If $\vv{g}\in (L^2)^m$ then
$M_{\vv{g}^\ast} $ is bounded on $\fM$ if and only if
\bequ\label{bisarasony}
\vv{g}\in (L^\infty)^m + \fM^\bot .
\eequ
Moreover, if $\vv{g} \in \fM^\bot$ then $M_{\vv{g}^\ast}\equiv 0$ on $\fM$.
\end{CMS-lemma}
\bdem
Let $\sigma\in (L^\infty)^{m\times m}$ be the matrix that satisfies $M_\sigma = P_{\fM}$, 
the orthogonal projection.
If $\vv{g}\in (L^2)^m$ is such that $M_{\vv{g}^\ast} $ is bounded on $\fM$ then
$M_{\vv{g}^\ast}  P_{\fM}  = M_{\vv{g}^\ast  \, \sigma}$ 
 is bounded on $(L^2)^m$, which yields $\vv{g}^\ast \, \sigma   \in (L^\infty)^{1\times m}$.
If $\sigma^\bot\in (L^\infty)^{m\times m}$ is such that $M_{\sigma^\bot}=P_{\fM^\bot}$ then
$$
\vv{g} = \sigma \vv{g} + \sigma^\bot \vv{g} =  (\vv{g}^\ast \sigma)^\ast + \sigma^\bot \vv{g}
\in (L^\infty)^{m\times 1} + \fM^\bot,
$$
as claimed. The other assertions follow from \eqref{bibess0}.
\edem

\noi
Suppose that $\fM\subset (L^2)^m$ is doubly invariant and $G \in  (L^2)^{m\times k}$ is a matrix
with columns $\vv{g}_j \in \fM$, for $1\le j\le k$.
By \eqref{bibess2-frame} the family
\bequ\label{bifamiglia}
\mathcal{F}:=\{ \S^n \vv{g}_j : \ n\in \bZ , \ 1\le j\le k \}
\eequ
is a Bessel sequence for $\fM$ if and only if  $M_{G^\ast} \big|_{\fM}$ is bounded,
and a frame if and only if it is  bounded above and below.

In  Theorem 2.3 of \cite{bow} Bownik characterized the families $\mathcal{F}$
that form a Bessel sequence, or a frame for its closed span  in terms of the fibers
$$\Phi(z)= \{ \vv{g}_j (z) : \  1\le j\le k,\ z\in \mathbb{T}\} \subset  \mathbb{C}^m.
$$
Namely, $\cF$ is a Bessel sequence with constant $B$ or  a frame with constants $A,\ B$
if and only if the same holds for $\Phi(z)$ uniformly for a.e $z\in \bT$.
Given a doubly invariant subspace $\fM \subset (L^2)^m$,
Theorem \ref{2coronaframe} below gives a characterization of the Bessel sequences and the frames
in $\fM$ in terms of the multiplication operator that projects on $\fM$.

\begin{CMS-propos}\label{bimanyy}
Let $\,\fM\subset (L^2)^m$ be a doubly invariant subspace and $\beta \in (L^\infty)^{m\times m}$
such that $M_\beta = P_{\fM^\bot}$.
If $F\in (L^\infty)^{m\times k}$, the following assertions are equivalent:
\begin{enumerate}
\item[{\em (1)}] $ M_{ FF ^\ast} \ge \delta^2 I_{(L^2)^m}$ on $\fM$ for some $\delta>0$.
\item[{\em (2)}] $F (z)  F ^\ast(z) + \beta(z) \ge \eta^2 I_{\mathbb{C}^m}$
for some $\eta>0$ and a.e.$\,z\in \bT$.
\end{enumerate}
\vspace{-2mm}
\end{CMS-propos}
\bdem
Since $\fM= \ker M_\beta$, it is clear that (1) is equivalent to
$M_{ FF ^\ast} + M_{\beta } \ge \varepsilon^2 I_{(L^2)^m}$ for some $\varepsilon>0$.
Also, this is an immediate consequence of (2).
Finally, if this inequality holds, consider a function $\chi_E(z) \vv{c}\in (L^2)^m$,
where $\vv{c}\in \mathbb{C}^m$ and $E\subset \bT$ is a set of measure $|E|>0$.
Then the above inequality says that
$$
\|\chi_E  F ^\ast \vv{c} \|^2_{(L^2)^k} + \|\chi_E  \beta \vv{c} \|^2_{(L^2)^m}
\ge \varepsilon^2 \|\chi_E\vv{c} \|^2_{(L^2)^m}= \varepsilon^2\| \vv{c} \|^2_{\mathbb{C}^m} \, |E|,$$
which can be rewritten as
$$
\int_E \big( \| F ^\ast(z) \vv{c} \|^2_{\mathbb{C}^k} + \|\beta(z) \vv{c} \|^2_{\mathbb{C}^m} \big)\, dz
\ge  \int_E  \varepsilon^2\| \vv{c} \|^2_{\mathbb{C}^m}\, dz,$$
where $dz$ denotes the Lebesgue measure on $\bT$. Since $E$ is arbitrary, we conclude that
(2) holds with $\eta = \varepsilon$.
\edem

\begin{CMS-theo}\label{2coronaframe}
Let $\fM\subset (L^2)^m$   be a doubly invariant subspace, $M_\sigma = P_{\fM}$
(with $\sigma\in (L^\infty)^{m\times m}$), and $g_j \in \fM$, where $j=1, \ldots , k$.
The family  $\cF$ in \eqref{bifamiglia} is a Bessel sequence  if and only if
$\vv{g}_j= M_\sigma \vv{f}_j$, with $\vv{f}_j\in (L^\infty)^m$ for  $1\le j\le k$,
and a frame for $\fM$ if and only if
$F:= \big[\vv{f}_{1}, \ldots , \vv{f}_{k}\big]\in  (L^\infty)^{m\times k}$
satisfies  Proposition \ref{bimanyy} with $\beta= I_{(L^2)^m}-\sigma$.
\end{CMS-theo}
\bdem
The assertion on Bessel sequences follows immediately from
\eqref{bibess2-frame} and \eqref{bisarasony}.
For the other assertion write $G:= \big[\vv{g}_{1}, \ldots , \vv{g}_{k}\big]$ and notice that
Lemma \ref{biDon} implies that $M_{G^\ast} \big|_{\fM} =M_{F^\ast} \big|_{\fM}$.
Thus, the comments following \eqref{bifamiglia} say that $\cF$ is a frame for $\fM$ if and only if
we have  (1) of Proposition \ref{bimanyy}.
\edem

\noi
In particular,  $\{U^nf, \ n\in \bZ\}$, with $f\in (L^2)^m$,
is a Bessel sequence if and only if $f\in (L^\infty)^m$.
The bilateral frame number of a model operator on $\fM\subset(L^2)^m$ can be obtained from
\cite[Thm.$\ $3.3]{bow} (see also the subsequent Remark (i) in\cite{bow}).
We give below a very simple proof that exploits a classical diagonalization result.

\begin{CMS-theo}[The bilateral frame number] \label{minimum-double}
Let $0\ne\fM\subset  (L^2)^m$ be a doubly invariant subspace and
$\Theta \in (L^\infty)^{m\times m}$ be  the matrix that satisfies $M_\Theta = P_{\fM^\bot}$.
Then
{\em
$$
\vartheta_- (\S |_\fM) =  \sup_{z\in \bT }  \big(\! \dim \ker \Theta(z)\big) ,\vspace{-2mm}
$$
}
where $\sup$ is the essential supremum and $\Theta(z)$ acts on $\bC^m$ for a.e.$\, z\in \bT$.
\end{CMS-theo}
\bdem
Since $\Theta$ is self-adjoint, a result of Kadison \cite[Thm.$\ $3.19]{Ka84} asserts that
there is a unitary matrix $V\in (L^\infty)^{m\times m}$ such that 
$V^\ast \Theta V = D \in (L^\infty)^{m\times m}$
is the diagonal matrix of a projection for a.e.$\, z\in \bT$.
Therefore, when $G\in (L^\infty)^{m\times k}$, we have
\begin{align*}
\inf_{ \vv{c}\in\bolam}  \| G^\ast (z) \vv{c} \|_{\bC^k}^2  + \| D (z) \vv{c} \|_{\bC^m}^2
&= \inf_{ \vv{c}\in\bolam}  \| G^\ast (z) V^\ast(z)\vv{c} \|_{\bC^k}^2  + \| D (z) V^\ast(z)\vv{c} \|_{\bC^m}^2 \\
&= \inf_{ \vv{c}\in\bolam}  \| \underbrace{G^\ast (z) V^\ast(z)}_{F^\ast(z)} \vv{c} \|_{\bC^k}^2  +
\| \underbrace{V(z) D (z) V^\ast(z)}_{\Theta(z)}   \vv{c} \|_{\bC^m}^2,
\end{align*}
where $F\in (L^\infty)^{m\times k}$.
Hence, if $\fD\subset (L^2)^m$ is the doubly invariant subspace that satisfies $M_D= P_{\fD^\bot}$,
we get $\vartheta_- (\S |_\fM) = \vartheta_-(\S|_\fD)$.
The diagonal of $D$ is formed by the characteristic functions of measurable sets
$E_1, \ldots , E_m \subset \bT$.
Since $G^\ast$ must satisfy
\bequ\label{maybelast}
\inf_{ \vv{c}\in\bolam}  \| G^\ast (z) \vv{c} \|_{\bC^k}^2  + \| D (z) \vv{c} \|_{\bC^m}^2 \ge \eta
\ \mbox{ for a.e.$\, z\in \bT$ and some $\eta>0$},
\eequ
the rank of $G^\ast (z)$ must be at least the cardinal of $\{ 1\le  j\le m:\ \chi_{E_j}(z) =0\}$. Thus,
$$
\vartheta_-(\S|_\fD) \ge \sup_{z\in \bT } \big(\! \dim \ker D(z) \big).\vspace{-1mm}
$$
For the other inequality let $p$ be the above supremum 
and let $1\le j_1 < \ldots < j_t\le m$ (with $t\le p$),
be such that $\chi_{E_{j_i}}(z) =0$.
Then the matrix $G^\ast \in (L^\infty)^{p\times m}$ whose rows are
$\vv{e}_{j_1}, \ldots , \vv{e}_{j_t}, 0 \ldots 0$, where $\vv{e}_j \in \bC^m$ is the standard  basis,
satisfies \eqref{maybelast}. Hence,
$$
\vartheta_- (\S |_\fM)= \vartheta_-(\S|_\fD) = \sup_{z\in \bT }   \big(\! \dim \ker D(z) \big)
= \sup_{z\in \bT }  \big(\! \dim \ker \Theta(z)\big).\vspace{-7mm}
$$
\edem

\subsection*{Acknowledgement}
 We are grateful to E. Andruchow for bringing Kadison's paper \cite{Ka84}, to our attention.

The research of the authors is partially supported by grants: UBACyT 20020170100430BA,  PICT 2018-3399 (ANPCyT)  and CONICET PIP 11220150100355.

%
%
%

\end{document}